\documentclass[11pt]{article}
\def\R#1#2#3#4{\hat{R}^{#1#2}_{#3#4}}
\def\Rm#1#2#3#4{\hat{R}^-{}^{#1#2}_{#3#4}}
\def\Rl#1#2#3#4{\grave{R}^{#1#2}_{#3#4}}
\def\Rlm#1#2#3#4{\grave{R}^-{}^{#1#2}_{#3#4}}
\def\Rr#1#2#3#4{\acute{R}^{#1#2}_{#3#4}}
\def\Rrm#1#2#3#4{\acute{R}^-{}^{#1#2}_{#3#4}}
\def\Rc#1#2#3#4{\check{R}^{#1#2}_{#3#4}}
\def\Rcm#1#2#3#4{\check{R}^-{}^{#1#2}_{#3#4}}

\def\d{\mathrm{d}}
\def\id{\mathrm{id}}

\def\nat{\mathrm{I\kern-.18em N}}
\def\real{\mathrm{I\kern-.18em R}}
\def\cplx{\mathrm{C\kern-0.5em\rule[0.2ex]{0.05em}{1.2ex}\kern0.6em}}
\def\SUq{\mathrm{SU}_q}
\def\Sq{\mathrm{S}_q}
\def\DeltaR{\Delta_\mathrm{R}}
\def\PhiR{\Phi_\mathrm{R}}
\def\s#1{\mathbf{s}_{#1}}
\def\se#1{{\mathbf{s}'_{#1}}}
\def\H#1{\Omega_{#1}}
\def\HtH{\H{}\otimes\H{}}
\def\sumL#1{\sum\limits_{#1=1}^N q^{-2#1}}
\def\sumR#1{\sum\limits_{#1=1}^N}
\newcommand{\Gammad}{{(\Gamma,\d)}}
\newcommand{\Gammat}{\tilde{\Gamma}}
\newcommand{\Gammatd}{{(\Gammat,\d)}}
\newcommand{\Gammaat}{\Gamma_{\alpha\tau}}
\newcommand{\Gammaatd}{{(\Gammaat,\d)}}
\newcommand{\Gammaap}{\Gamma'_{\alpha\omega}}
\newcommand{\Gammaapd}{{(\Gammaap,\d)}}
\newcommand{\Gammap}{\Gamma^{\prime\prime}_{\omega\psi}}
\newcommand{\Gammapd}{{(\Gammap,\d)}}
\newcommand{\Gammapp}{\Gamma^{\prime\prime\prime}_{\varrho\tau}}
\newcommand{\Gammappd}{{(\Gammapp,\d)}}
\newcommand{\Gammatl}[1][\lambda]{\tilde{\Gamma}_{#1}}
\newcommand{\Gammatld}[1][\lambda]{{(\Gammatl[#1],\d)}}
\newcommand{\Gammatp}[1][\lambda]{\tilde{\Gamma}^{\prime}_{#1}}
\newcommand{\Gammatpd}[1][\lambda]{{(\Gammatp[#1],\d)}}
\newcommand{\Gammatpp}[1][\lambda]{\tilde{\Gamma}^{\prime\prime}_{#1}}
\newcommand{\Gammatppd}[1][\lambda]{{(\Gammatpp[#1],\d)}}
\newcommand{\Gammato}[1][\lambda]{\tilde{\Gamma}^{\bullet}_{#1}}
\newcommand{\Gammatod}[1][\lambda]{{(\Gammato[#1],\d)}}
\newcommand{\Gammatoo}[1][\lambda]{\tilde{\Gamma}^{\bullet\bullet}_{#1}}
\newcommand{\Gammatood}[1][\lambda]{{(\Gammatoo[#1],\d)}}
\newtheorem{Prop}{Proposition}
\newtheorem{Thm}[Prop]{Theorem}
\newtheorem{Lemma}[Prop]{Lemma}
\newtheorem{Cor}[Prop]{Corollary}
\newcounter{DefC}
\newenvironment{Def}
    {\stepcounter{DefC}\vskip1ex\textbf{Definition \arabic{DefC}~}}
    {\par\vspace{1.5ex}}
\def\qed{\hbox to1cm{~}\hspace*\fill\rule{1ex}{1ex}\par\vspace{1ex}}
\newenvironment{Proof}
    {\textbf{Proof}~}
    {\qed}

\def\cc#1{\kern0.1em%
  \overline{\kern-0.1em#1\rule{0pt}{1.5ex}\kern-0.1em}\kern0.1em}
\def\calpha{{\cc{\alpha}}}
\def\rot{{\displaystyle\frac{\varrho}{\tau}}}
\def\tor{{\displaystyle\frac{\tau}{\varrho}}}
\begin{document}
\centerline{\LARGE\bf Differential Calculus on Quantum Spheres}

\begin{center}
Martin Welk\\[2ex]
Universit{\"a}t Leipzig, Institut f{\"u}r Mathematik\\
Augustusplatz 9, 04109 Leipzig\\Germany
\end{center}
\begin{abstract}
We study covariant differential calculus on the quantum spheres $\Sq^{2N-1}$.
Two classification results for covariant first order differential calculi
are proved. As an important step towards a description of the noncommutative
geometry of the quantum spheres, a framework of covariant differential
calculus is established, including a particular first order calculus obtained
by factorisation, higher order calculi and a symmetry concept.
\end{abstract}

\section{Introduction}
Quantum groups and quantum spaces are
important examples of noncommutative geometric spaces.
The description of differential calculus on them
forms the fundament for an analysis of their geometric structure.

Quantum groups are the most advanced object of study.
Covariant---and especially, bicovariant---differential calculi on
quantum groups have been under investigation during the last years
\cite{wr,ss1,ss2},
and basic concepts of differential geometry on quantum groups have already
been introduced, see e.\,g.\ \cite{iv}.
There are also results concerning covariant differential calculi on
several examples of quantum spaces, e.\,g.\ quantum vector spaces,
and Podle{\'s}' spheres \cite{ap}.

Quantum homogeneous spaces are a class of quantum spaces which are
in an especially close relation to quantum groups, therefore presenting
themselves as a promising object for investigation.

In this paper, we study the quantum spheres $\Sq^{2N-1}$
as introduced by Vaksman and Soibelman \cite{vs} as an example of
a quantum homogeneous space.

We start by studying covariant first order differential calculi on these
quantum spaces. We prove, as our main results on this topic,
two classification theorems for first order calculi under slightly
different selective constraints. The classification results hold for
$N\ge4$ but the differential calculi included exist for $N=2$ and $N=3$, too.
We point out the relations between the two sets of calculi.

Subsequently, we describe higher order differential calculus on the
quantum spheres based on a particular first order differential calculus.
Our approach uses ideas from the well-developed theory of bicovariant
differential calculi on quantum groups, thereby providing us even with
a symmetry concept for tensor products of differential forms linked
with the higher order calculus.
The framework of higher order differential calculus and symmetry
is powerful enough to enable the introduction of basic
concepts of noncommutative differential geometry on
the quantum spheres $\Sq^{2N-1}$.

\section{Quantum spaces}

\subsection{General definitions and conventions}
We start by collecting some basic definitions on quantum spaces,
mostly following the terminology as e.\,g.\ in~\cite{ap}.

Let ${\cal A}$ be a Hopf algebra with comultiplication $\Delta$ and
counit $\varepsilon$. A \emph{quantum space}
for ${\cal A}$ is a
pair $(X, \DeltaR)$ where $X$ is a unital algebra and
$\DeltaR: X\to X\otimes{\cal A}$ a \emph{(right) coaction of
${\cal A}$ on $X$}, i.\,e.\ an
algebra homomorphism such that $(\DeltaR\otimes\id)\DeltaR
=(\id\otimes\Delta)\DeltaR$; $(\id\otimes\varepsilon)\DeltaR
=\id$.
$X$ is called {\em quantum homogeneous space} for ${\cal A}$
if there is an embedding $\iota:X\to{\cal A}$ such that
$\DeltaR=\Delta\circ\iota$.

\vspace{1ex}

Throughout this paper, the dimension $N$ of the underlying quantum group
$\SUq(N)$ is a natural number, $N\ge2$. The deformation parameter
$q$ is a real number, ${q\not\in\{-1,0,1\}}$. We use the abbreviations
${Q=q-q^{-1}}$, ${\s+=\sum_{i=0}^{N-1}q^{2i}}$, ${\se+=\s+-1}$.

We also need the R-matrices which are well-known from
investigations on the quantum group $\SUq(N)$, e.\,g.\ \cite{rtf,ss1}.
Note that $\R{}{}{}{}$ is an invertible ${N^2\times N^2}$ matrix with
the inverse $\Rm{}{}{}{}$, and that ${\R{}{}{}{}-\Rm{}{}{}{}=Q I}$
with $I$ being the ${N^2\times N^2}$ unit matrix.
\[
\R ijkl = \left\{\begin{array}{l@{~~}l}
1&\mbox{for $i=l\ne k=j$}\\
q&\mbox{for $i=j=k=l$}\\
Q&\mbox{for $i=k<j=l$}\\
0&\mbox{otherwise}
\end{array}\right.
;\qquad
\Rm ijkl = \left\{\begin{array}{l@{~~}l}
1&\mbox{for $i=l\ne k=j$}\\
q^{-1}&\mbox{for $i=j=k=l$}\\
-Q&\mbox{for $i=k>j=l$}\\
0&\mbox{otherwise}
\end{array}\right.
\]

The following matrices are derived from these fundamental R-matrices:
\[
\begin{array}{@{}r@{~}l@{\quad}r@{~}l@{\quad}r@{~}l}
\Rc ijkl&= \R lkji;&
\Rl ijkl&= q^{2l-2i}\R jlik;&
\Rr ijkl&= \R kilj;\\[.5ex]
\Rcm ijkl&= \Rm lkji;&
\Rlm ijkl&= q^{2l-2i}\Rm jlik;&
\Rrm ijkl&= \Rm kilj.
\end{array}
\]

\subsection{The quantum spheres $\Sq^{2N-1}$}

Our object of study are the quantum spheres introduced
by Vaksman and Soibelman \cite{vs} which we shall describe now.

Let $X$ be the free unital algebra with a set of $2N$ generators
${\{z_i,z^*_i~|~ i=1,\dots,N\}}$ and defining relations
\begin{equation}\label{qsr1}
\begin{array}{rcl@{\quad}l}
z_iz_j&=&qz_jz_i&(1\le i<j\le N)\\[.5ex]
z^*_i z^*_j&=&q^{-1} z^*_j z^*_i&(1\le i<j\le N)\\[.5ex]
z_i z^*_j&=&q z^*_j z_i&(1\le i,j\le N,~i\ne j)
\end{array}
\end{equation}
\begin{equation}\label{qsr23}
z_i z^*_i- z^*_i z_i+q^{-1}Q\sum\limits_{k>i}z_k z^*_k=0,
\qquad\qquad
\sumR{i}z_i z^*_i=1.
\end{equation}

This algebra is made into a $*$~algebra by letting
$(z_i)^*=z^*_i$; $(z^*_i)^*=z_i$.
Then, $X$ is called \emph{quantum sphere} and denoted by $\Sq^{2N-1}$.

Using the R-matrices, the relations (\ref{qsr1}) can be rewritten as
\[
\begin{array}{rcl@{\qquad}rcl}
\R klij z_kz_l&=&qz_iz_j;&
\Rcm klij  z^*_k z^*_l&=&q^{-1} z^*_i z^*_j;\\[.5ex]
\Rr klij z_k z^*_l&=&q z^*_i z_j;&
\Rlm klij  z^*_k z_l&=&q^{-1}z_i z^*_j.
\end{array}
\]

The relations (\ref{qsr23}) imply
\[
\sumL{i}z^*_i z_i=q^{-2}\quad\textrm{and}\quad
z_iz^*_i-z^*_iz_i+qQ\sum\limits_{k>i}q^{2i-2k}z^*_kz_k=0.\]

Let $u^i_j$, $1\le i,j\le N$
be the generators and $S$ the antipode map of the quantum group
$\SUq(N)$ as defined in \cite{rtf}.
Then, by
\[z_i=u^1_i,\qquad z^*_i=(u^1_i)^*=S(u^i_1)\]
an embedding of $\Sq^{2N-1}$ into $\SUq(N)$ is given, making the quantum
sphere into a quantum homogeneous space for $\SUq(N)$ with the coaction
\[\DeltaR(z_i)=\sumR{j}z_j\otimes u^j_i;\qquad
  \DeltaR(z^*_i)=\sumR{j}z^*_j\otimes S(u^i_j).\]

\section{First order differential calculus}
\subsection{Basic definitions}
First we recall important definitions concerning first order differential
calculi, cf.~\cite{ap}.

A \emph{first order differential calculus} on an algebra $X$ means
a pair $(\Gamma, \d)$ where $\Gamma$ is a bimodule over $X$ and
$\d:X\to\Gamma$ is a linear mapping which fulfils
Leibniz' rule $\d(xy)=(\d x)y+x(\d y)$ for all
$x,y\in X$, and $\Gamma=\mathrm{Lin}\{x\d y~ |~ x,y\in X\}$.
The elements of $\Gamma$ are called one-forms.

A first order differential calculus $(\Gamma,\d)$
on a quantum space $(X,\DeltaR)$ for ${\cal A}$ is called
\emph{(right) covariant} if there is a linear mapping
$\PhiR:\Gamma\to\Gamma\otimes{\cal A}$
with $(\PhiR\otimes\id)\PhiR=(\id\otimes\Delta)\PhiR$;
$(\id\otimes\varepsilon) \PhiR=\id$;
$\PhiR(x\omega y)=\DeltaR(x)\PhiR(\omega)\DeltaR(y)$; $\PhiR(\d x)=(\d\otimes
\id)\DeltaR(x)$ for all $x,y\in X$, $\omega\in\Gamma$.
A one-form $\omega$ is called \emph{invariant} if
${\DeltaR(\omega)=\omega\otimes1}$.

If $X$ is a $*$~algebra then
a first order differential calculus $(\Gamma, \d)$ is called a
\emph{$*$~calculus}
if $\sum\limits_k\,x_k\d y_k=0$ implies
$\sum\limits_k\,\d(y_k^*)x_k^*=0$ for $x_k,y_k\in X$.

\subsection{First order differential calculi on quantum spheres---Results}
\label{srs1}

In this section, we give two classification results for covariant
first order differential calculi on the quantum spheres $\Sq^{2N-1}$
differing by the set of selective constraints used for
classification.

Our first theorem gives a classification of covariant first order
differential $*$~calculi on the quantum spheres with the
constraint that the calculi be freely generated as left modules
by the differentials of the generators of the quantum sphere.

Most of these calculi allow for a factorisation by an additional relation
${\H0=0}$ where $\H0$ is an invariant one-form of $\Gamma$, yielding
first order differential $*$~calculi of a different kind.
We shall see that a relation of this kind holds in the classical case, too.
Therefore it is to be expected that differential calculi of this
second type are more adequate
to describe the noncommutative geometry of the quantum spheres than
are the freely generated ones. This leads us to give also a direct
classification of first order differential calculi of this second kind.
We shall relax the selective constraint of our classification in this
second case to give full account not only of $*$~calculi but of all
covariant first order differential calculi on $\Sq^{2N-1}$ for which
all relations in the left module $\Gamma$ are algebraically generated by
one relation ${\H0=0}$ where $\H0$ is a fixed invariant one-form.

It is clear that in any covariant first order differential
calculus~$\Gammad$ on $\Sq^{2N-1}$ the two one-forms
\begin{equation}
\H+=\sumR{i}z_i\d z^*_i,\qquad\H-=\sumL{i}z^*_i\d z_i
\end{equation}
are invariant and that any invariant one-form in $\Gamma$ is a linear
combination of $\H+$ and $\H-$.

\begin{Thm}\label{fodc}
On $\Sq^{2N-1}$, there exist first order differential $*$~calculi
\begin{itemize}
\item ${\Gammad=\Gammaatd}$ where ${\alpha\in\real\setminus\{0,q^{-2}\}}$ and
${\tau\in\real}$,
\item ${\Gammad=\Gammaapd}$
where either ${\alpha\in\real\setminus\{0,q^{-2}\}}$ and
${\omega\in\real\setminus\{0\}}$, or ${\alpha\in\cplx\setminus\real}$ and
${\omega=q^4\alpha\calpha}$,
\item ${\Gammad=\Gammapd}$ where ${\omega\in\real\setminus\{0\}}$ and
${\psi\in\real}$,
\item ${\Gammad=\Gammappd}$ where ${\varrho,\tau\in\real\setminus\{0\}}$,
\end{itemize}
which are covariant with respect to $\SUq(N)$ and for which
${\{\d z_i,\d z^*_i~|~i=1,\dots,N\}}$ is a free left module basis for
$\Gamma$, with their bimodule structure given by
\begin{eqnarray*}
\hbox to2cm{$\Gammaat:$\hfil
$\d z_k z_l$}&=&q\alpha\Rm stkl z_s\d z_t+(q^2\alpha-1)z_k\d z_l\\*
  &&\qquad+q^2\alpha^2(1-\se+\tau)z_kz_l\H+
  +q^2(1-\alpha\se+\tau)z_kz_l\H-\\
\d z^*_k z^*_l&=&q^{-1}\alpha^{-1}\Rc stkl z^*_s\d z^*_t
  +(q^{-2}\alpha^{-1}-1)z^*_k\d z^*_l\\*
  &&\qquad+(1-\se+\tau)z^*_kz^*_l\H+
  +\alpha^{-2}(1-\alpha\se+\tau)z^*_kz^*_l\H-\\
\d z_k z^*_l&=&q^{-1}\alpha^{-1}\Rlm stkl z^*_s\d z_t
  +(q^2\alpha-1)z_k\d z^*_l-q^2\alpha(1-\s+\tau)z_kz^*_l\H+\\*
  &&\qquad-\alpha\tau q^{2k}\delta_{kl}\H+
  -\alpha^{-1}(1-q^2\alpha\s+\tau)z_kz^*_l\H--\tau q^{2k}\delta_{kl}\H-\\
\d z^*_k z_l&=&q\alpha\Rr stkl z_s\d z^*_t
  +(q^{-2}\alpha^{-1}-1)z^*_k\d z_l-q^2\alpha(1-\s+\tau)z^*_kz_l\H+\\*
  &&\qquad-\alpha\tau\delta_{kl}\H+
  -\alpha^{-1}(1-q^2\alpha\s+\tau)z^*_kz_l\H--\tau\delta_{kl}\H-\\[.7ex]
\hbox to2cm{$\Gammaap:$\hfil
$\d z_k z_l$}&=&q\alpha\Rm stkl z_s\d z_t+(q^2\alpha-1)z_k\d z_l\\*
&  &\qquad+\omega z_kz_l\H+
  +(\alpha^{-1}\omega-q^2(\alpha-1))z_kz_l\H-\\
\d z^*_k z^*_l&=&q^{-1}\calpha^{-1}\Rc stkl z^*_s\d z^*_t
  +(q^{-2}\calpha^{-1}-1)z^*_k\d z^*_l\\*
&  &\qquad+(q^2\calpha\omega^{-1}-(\calpha^{-1}-1))z^*_kz^*_l\H+
  +q^2\omega^{-1}z^*_kz^*_l\H-\\
\d z_k z^*_l&=&q^{-1}\alpha^{-1}\Rlm stkl z^*_s\d z_t
  +(q^2\calpha-1)z_k\d z^*_l\\*
&  &\qquad-q^2\calpha z_kz^*_l\H+-\alpha^{-1}z_kz^*_l\H-\\
\d z^*_k z_l&=&q\calpha\Rr stkl z_s\d z^*_t
  +(q^{-2}\alpha^{-1}-1)z^*_k\d z_l\\*
&  &\qquad-q^2\calpha z^*_kz_l\H+-\alpha^{-1}z^*_kz_l\H-\\[.7ex]
\hbox to2cm{$\Gammap:$\hfil
$\d z_k z_l$}&=&q^{-1}\Rm stkl z_s\d z_t+\omega z_kz_l\H+
  +(q^2\omega\psi-1)z_kz_l\H-\\
\d z^*_k z^*_l&=&q\Rc stkl z^*_s\d z^*_t+(\psi-q^2)z^*_k z^*_l\H+
  +q^2\omega^{-1}z^*_kz^*_l\H-\\
\d z_k z^*_l&=&q\Rlm stkl z^*_s\d z_t-z_kz^*_l\H+-q^2z_kz^*_l\H-\\
\d z^*_k z_l&=&q^{-1}\Rr stkl z_s\d z^*_t-z^*_k z_l\H+-q^2z^*_kz_l\H-\\[.7ex]
\hbox to2cm{$\Gammapp:$\hfil
$\d z_k z_l$}&=&q^{-1}\Rm stkl z_s\d z_t
  -q^{-2}\rot(\se+\varrho-1)z_kz_l\H+-\rot(\se+\tau-q^2)z_kz_l\H-\\
\d z^*_k z^*_l&=&q\Rc stkl z^*_s\d z^*_t
  -\tor(\se+\varrho-1)z^*_kz^*_l\H+-q^2\tor(\se+\tau-q^2)z^*_kz^*_l\H-\\
\d z_k z^*_l&=&q\Rlm stkl z^*_s\d z_t
  -q^{-2}\varrho q^{2k}\delta_{kl}\H+-\tau q^{2k}\delta_{kl}\H-\\*
&  &\qquad+(\s+\varrho-1)z_kz^*_l\H++q^2(\s+\tau-1)z_kz^*_l\H-\\
\d z^*_k z_l&=&q^{-1}\Rr stkl z_s\d z^*_t
  -q^{2N-2}\varrho\delta_{kl}\H+-q^{2N}\tau\delta_{kl}\H-\\*
&  &\qquad+(\s+\varrho-1)z^*_kz_l\H++q^2(\s+\tau-1)z^*_kz_l\H-.
\end{eqnarray*}
If $N\ge4$, any first order differential $*$~calculus $\Gammad$
on $\Sq^{2N-1}$ which is covariant with respect to $\SUq(N)$ and for which
$\{\d z_i,\d z^*_i~|~ i=1,\dots,N\}$ is a free left module basis for
$\Gamma$ is isomorphic to one of the calculi $\Gammaatd$, $\Gammaapd$,
$\Gammapd$, $\Gammappd$.
\end{Thm}

None of the first order differential calculi from the theorem is inner
since there exists
no invariant one-form $\tilde{\H{}}=\alpha\H++\beta\H-$
in $\Gamma_{\omega}$ for $\omega\in\real\setminus\{0\}$ such that
$\d x=\tilde{\H{}}x-x\tilde{\H{}}$ for all $x\in\Sq^{2N-1}$.
Instead, one easily checks that in all calculi, except for $\Gammaapd$ with
non-real $\alpha$, there is at least one
invariant element ${\H0=\mu\H++\lambda\H-}$ which quasi-commutes
with any $x\in\Sq^{2N-1}$, i.\,e.\ ${\H0x=C(x)x\H0}$
where $C(x)$ is a complex number dependent on $x\in X$.

\begin{Cor}\label{inv}
In $\Gammaatd$, and $\Gammaapd$ with real $\alpha$,
there exists (up to scalar multiples) exactly one invariant one-form
${\H++\alpha^{-1}\H-}$ that quasi-commutes with all ${x\in X}$.

In $\Gammapd$, there exist (up to scalar multiples)
exactly two invariant one-forms that quasi-commute with all ${x\in X}$.
They are given by ${\H++\lambda_1\H-}$ and ${\H++\lambda_2\H-}$
where $\lambda_{1,2}$ are the solutions of the quadratic equation
\[q^{-4}\omega\lambda^2-q^{-2}\omega\psi+1=0.\]

In $\Gammappd$, with ${\tau/\varrho\ne q^2}$,
there exist (up to scalar multiples)
exactly two invariant one-forms that quasi-commute with all ${x\in X}$.
One of them is given by ${\H++\lambda_1\H-}$, ${\lambda_1=q^2\tau/\varrho}$,
while the other one is given by $\H+$, if ${\se+\varrho=1}$, or by
${\H++\lambda_2\H-}$,
${\lambda_2=q^2(\se+\tau-q^2)/(\se+\varrho-1)}$, otherwise.

In $\Gammappd$, with ${\varrho=q^{-2}\tau\ne\se+^{-1}}$,
there exists (up to scalar multiples) exactly one invariant one-form,
${\H++q^4\H-}$, that quasi-commutes with all ${x\in X}$.

In $\Gammappd$, with ${\varrho=q^{-2}\tau=\se+^{-1}}$, any invariant one-form
quasi-commutes with all ${x\in X}$.
\end{Cor}
\begin{Proof}
By direct calculation, one easily checks the quasi-commutation
statements for the given invariant one-forms with $z^{(*)}_m$.

By transforming ${(\mu\H++\lambda\H-)z^{(*)}_m}$, with variable
coefficients $\mu$, $\lambda$, into a left module expression,
we derive necessary conditions for $\mu$ and $\lambda$ which lead
to the uniqueness assertions of the Corollary.
\end{Proof}

If, for a first order differential calculus $\Gammad$, the invariant
one-form $\H0$ quasi-commutes with all $x\in X$, the calculus
$\Gammad$ allows for a factorisation by the additional relation ${\H0=0}$.
Therefore, we obtain from the calculi of Theorem~\ref{fodc} new calculi
$\Gammatd$ for which ${\{\d z_i,\d z^*_i~|~i=1,\dots,N\}}$ is no longer a free
left module basis for~$\Gammat$.

A factorisation of this kind is quite natural since
in the classical limit (${q=1}$) we have
$\H++\H-=\d\left(\sumR{i} z_iz^*_i\right)=0$ anyway.

\begin{Thm}\label{fodcf}
On $\Sq^{2N-1}$, there exist first order differential calculi
\begin{itemize}
\item $\Gammatd=\Gammatld$ where ${\lambda\in\cplx\setminus\{0,q^2\}}$,
\item $\Gammatd=\Gammatpd$ where ${\lambda\in\cplx\setminus\{0\}}$,
\item $\Gammatd=\Gammatppd$ where ${\lambda\in\cplx\cup\{\infty\}}$,
\item $\Gammatd=\Gammatod$ where ${\lambda\in\cplx\setminus\{0,q^2\}}$,
\item $\Gammatd=\Gammatood$ where ${\lambda\in\cplx\setminus\{0,q^2\}}$,
\end{itemize}
which are covariant with respect to $\SUq(N)$
and for which all relations
in the left module $\Gammat$ are algebraically generated by
one relation ${\H++\lambda\H-=0}$ (if ${\lambda\in\cplx}$)
or ${\H-=0}$ (if ${\lambda=\infty}$),
with the bimodule structure given by
\begin{eqnarray*}
\hbox to2cm{$\Gammatl:$\hfil
$\d z_kz_l$}&=&q\lambda^{-1}\Rm stkl z_s\d z_t+(q^2\lambda^{-1}-1)z_k\d z_l
   +q^2\lambda^{-1}(\lambda^{-1}-1)z_kz_l\H+\\
\d z^*_kz^*_l&=&q^{-1}\lambda\Rc stkl z^*_s\d z^*_t
   +(q^{-2}\lambda-1)z^*_k\d z^*_l-(\lambda-1)z^*_kz^*_l\H+\\
\d z_kz^*_l&=&q^{-1}\lambda\Rlm stkl z^*_s\d z_t
   +(q^2\lambda^{-1}-1)z_k\d z^*_l-(q^2\lambda^{-1}-1)z_kz^*_l\H+\\
\d z^*_kz_l&=&q\lambda^{-1}\Rr stkl z_s\d z^*_t
   +(q^{-2}\lambda-1)z^*_k\d z_l-(q^2\lambda^{-1}-1)z^*_kz_l\H+\\[.7ex]
\hbox to2cm{$\Gammatp:$\hfil
$\d z_kz_l$}&=&q^{-1}\Rm stkl z_s\d z_t
   -\lambda^{-1}(q^4\lambda^{-1}-1)z_kz_l\H+\\
\d z^*_kz^*_l&=&q\Rc stkl z^*_s\d z^*_t
   -q^{-2}\lambda(q^4\lambda^{-1}-1)z^*_kz^*_l\H+\\
\d z_kz^*_l&=&q\Rlm stkl z^*_s\d z_t
   +(q^2\lambda^{-1}-1)z_kz^*_l\H+\\
\d z^*_kz_l&=&q^{-1}\Rr stkl z_s\d z^*_t
   +(q^2\lambda^{-1}-1)z^*_kz_l\H+\\[.7ex]
\hbox to2cm{$\Gammatpp,~\lambda\not\in\{0,\infty\}:$\kern-2cm\hfil}\\*
\d z_kz_l&=&q^{-1}\Rm stkl z_s\d z_t\\
\d z^*_kz^*_l&=&q\Rc stkl z^*_s\d z^*_t\\
\d z_kz^*_l&=&q\Rlm stkl z^*_s\d z_t
   +q^{-2}\se+^{-1}(q^4\lambda^{-1}-1)q^{2k}\delta_{kl}\H+\\*&&\qquad
   -\se+^{-1}(q^{2N+2}\lambda^{-1}-1)z_kz^*_l\H+\\
\d z^*_kz_l&=&q^{-1}\Rr stkl z_s\d z^*_t
   +q^{-2}\se+^{-1}(q^4\lambda^{-1}-1)\delta_{kl}\H+\\*&&\qquad
   -\se+^{-1}(q^{2N+2}\lambda^{-1}-1)z^*_kz_l\H+\\[.7ex]
\hbox to2cm{$\Gammatpp[0]:$\hfil
$\d z_kz_l$}&=&q^{-1}\Rm stkl z_s\d z_t\\
\d z^*_kz^*_l&=&q\Rc stkl z^*_s\d z^*_t\\
\d z_kz^*_l&=&q\Rlm stkl z^*_s\d z_t
   -q^{-2N+2}\se+^{-1}q^{2k}\delta_{kl}\H-
   +q^2\se+^{-1}z_kz^*_l\H-\\
\d z^*_kz_l&=&q^{-1}\Rr stkl z_s\d z^*_t
   -q^2\se+^{-1}\delta_{kl}\H-
   +q^2\se+^{-1}z^*_kz_l\H-\\[.7ex]
\hbox to2cm{$\Gammatpp[\infty]:$\hfil
$\d z_kz_l$}&=&q^{-1}\Rm stkl z_s\d z_t\\
\d z^*_kz^*_l&=&q\Rc stkl z^*_s\d z^*_t\\
\d z_kz^*_l&=&q\Rlm stkl z^*_s\d z_t
   -q^{-2}\se+^{-1}q^{2k}\delta_{kl}\H+
   +\se+^{-1}z_kz^*_l\H+\\
\d z^*_kz_l&=&q^{-1}\Rr stkl z_s\d z^*_t
   -q^{2N-2}\se+^{-1}\delta_{kl}\H+
   +\se+^{-1}z^*_kz_l\H+\\[.7ex]
\hbox to2cm{$\Gammato:$\hfil
$\d z_kz_l$}&=&q\lambda^{-1}\Rm stkl z_s\d z_t+(q^2\lambda^{-1}-1)z_k\d z_l
   +q^2\lambda^{-1}(\lambda^{-1}-1)z_kz_l\H+\\
\d z^*_kz^*_l&=&q\Rc stkl z^*_s\d z^*_t
   -q^{-2}\lambda(q^4\lambda^{-1}-1)z^*_kz^*_l\H+\\
\d z_kz^*_l&=&q^{-1}\lambda\Rlm stkl z^*_s\d z_t\\
\d z^*_kz_l&=&q^{-1}\Rr stkl z_s\d z^*_t
   +(q^{-2}\lambda-1)z^*_k\d z_l\\[.7ex]
\hbox to2cm{$\Gammatoo:$\hfil
$\d z_kz_l$}&=&q^{-1}\Rm stkl z_s\d z_t
   -\lambda^{-1}(q^4\lambda^{-1}-1)z_kz_l\H+\\
\d z^*_kz^*_l&=&q^{-1}\lambda\Rc stkl z^*_s\d z^*_t
   +(q^{-2}\lambda-1)z^*_k\d z^*_l-(\lambda-1)z^*_kz^*_l\H+\\
\d z_kz^*_l&=&q\Rlm stkl z^*_s\d z_t
   +(q^2\lambda^{-1}-1)z_k\d z^*_l\\
\d z^*_kz_l&=&q\lambda^{-1}\Rr stkl z_s\d z^*_t.
\end{eqnarray*}

If $N\ge4$, any first order differential calculus $\Gammatd$
on $\Sq^{2N-1}$ which is covariant with respect to $\SUq(N)$
and for which all relations
in the left module $\Gammat$ are algebraically generated by
one relation ${\H0=0}$ with invariant ${\H0\in\Gammat}$
is isomorphic to one of the calculi
$\Gammatld$, $\Gammatpd$, $\Gammatppd$, $\Gammatod$, $\Gammatood$.

The differential calculi
$\Gammatld$ for ${\lambda\in\real\setminus\{0,q^2\}}$,
$\Gammatpd$ for ${\lambda\in\real\setminus\{0\}}$,
and $\Gammatppd$ for ${\lambda\in\real\cup\{\infty\}}$ are $*$~calculi.

If $N\ge4$, any first order differential $*$~calculus $\Gammatd$
on $\Sq^{2N-1}$ which is covariant with respect to $\SUq(N)$
and for which all relations
in the left module $\Gammat$ are algebraically generated by
one relation ${\H0=0}$ (${\H0\in\Gammat}$ invariant),
is isomorphic to one of these calculi.
\end{Thm}

\begin{Cor}\label{coi}
The calculi $\Gammatld$ are inner, with
\[\d x=\H{}x-x\H{}~~\mbox{for all $x\in X$}\]
for ${\H{}:=(q^{-2}\lambda-1)^{-1}}$.
None of the other calculi described in Theorem~\ref{fodcf} is inner.
\end{Cor}
\begin{Proof}
In each of the calculi of the Theorem, all invariant one-forms are
scalar multiples of one invariant element $\H+$ or $\H-$.
The statement for $\Gammatld$ is proved by calculation.
The assertion for the remaining calculi follows from the fact that
in each of them the non-zero invariant one-forms (which are scalar
multiples of one element, $\H+$ or $\H-$) quasi-commute with
the algebra generators from at least one of the sets
${\{z_i~|~i=1,\dots,N\}}$ or ${\{z^*_i~|~i=1,\dots,N\}}$.
\end{Proof}

Now we describe the factorisation of the freely generated
first order differential $*$~calculi from Theorem~\ref{fodc} by the
relations ${\H0=0}$ where $\H0$ are the quasi-commuting invariant one-forms
according to Corollary~\ref{inv}.

\begin{Cor}
From the first order differential $*$~calculi from Theorem~\ref{fodc},
the following first order differential $*$~calculi are obtained by
factorisation:
\begin{itemize}
\item for any ${\alpha\in\real\setminus\{0,q^{-2}\}}$,
\[
\Gammaat/(\H++\alpha^{-1}\H-)=
\Gammaap/(\H++\alpha^{-1}\H-)=
\Gammatl[\alpha^{-1}];
\]
\item for any ${\omega\in\real\setminus\{0\}}$, ${\psi\in\real}$, and
$\lambda_{1,2}$ as in Corollary~\ref{inv},
\[
\Gammap/(\H++\lambda_1\H-)=\Gammatp[\lambda_1]
~~\mbox{and}~~
\Gammap/(\H++\lambda_2\H-)=\Gammatp[\lambda_2];
\]
\item for any ${\varrho\in\real\setminus\{0,\se+^{-1}\}}$,
${\tau\in\real\setminus\{0\}}$, and $\lambda_{1,2}$ as in Corollary~\ref{inv}
(${\lambda_1=\lambda_2=q^4}$ if ${\tau=q^2\varrho}$),
\[
\Gammapp/(\H++\lambda_1\H-)=\Gammatp[\lambda_1]
~~\mbox{and}~~
\Gammapp/(\H++\lambda_2\H-)=\Gammatpp[\lambda_2];
\]
\item for ${\varrho=\se+^{-1}}$ and any
${\tau\in\real\setminus\{0,q^2\se+^{-1}\}}$,
and $\lambda_1$ as in Corollary~\ref{inv},
\[
\Gammapp/(\H++\lambda_1\H-)=\Gammatp[\lambda_1]
~~\mbox{and}~~
\Gammapp/(\H-)=\Gammatpp[\infty];
\]
\item for ${\varrho=q^{-2}\tau=\se+^{-1}}$, and any ${\lambda\in\real}$,
\[
\Gammapp/(\H++\lambda\H-)=\Gammatpp[\lambda];~~
\Gammapp/(\H-)=\Gammatpp[\infty].
\]
\end{itemize}
\end{Cor}
\begin{Proof}
By eliminating $\H-$ ($\H+$ for $\Gammapp$ with ${\lambda_2=0}$) from
the equations describing the bimodule structure of $\Gammaat$, $\Gammaap$,
$\Gammap$, $\Gammapp$, one obtains the corresponding equations for
the factorised calculi.
\end{Proof}

Let us now characterise one particular $*$~calculus which is of
particular importance for our further considerations.
It is the calculus $\Gammatl[1]$ with the relation
\begin{equation}
\H++\H-=\sumR{i}z_i\d z^*_i+\sumL{i}z^*_i\d z_i=0.
\label{fr}
\end{equation}
The bimodule structure of $\Gammatl[1]$ is given by
\begin{equation}
\begin{array}{@{}r@{~}l@{}}
\d z_kz_l&=q\R stkl z_s\d z_t\\[0.7ex]
\d z^*_kz^*_l&=q^{-1}\Rcm stkl z^*_s\d z^*_t\\[0.7ex]
\d z_kz^*_l&=q^{-1}\Rlm stkl z^*_s\d z_t
   +qQz_k\d z^*_l+Q^2z_kz^*_l\H{}\\[0.7ex]
\d z^*_kz_l&=q\Rr stkl z_s\d z^*_t
   -q^{-1}Qz^*_k\d z_l+Q^2z^*_kz_l\H{},
\end{array}
\label{fodcfb}
\end{equation}
where ${\H{}=-qQ^{-1}\H+=qQ^{-1}\H-}$.

Our considerations concerning higher order differential calculus and
symmetry in section~\ref{hos} will be based on this calculus
because it has two essential properties which are fulfilled simultaneously
only by this calculus.
\begin{itemize}
\item
First, $\Gammatl[1]$ is inner, as stated in Corollary~\ref{coi}. We have
\begin{equation}
\d x=\H{}x-x\H{}~~\mbox{for all $x\in X$.}
\label{foi}
\end{equation}
\item
Second, this calculus
decomposes into subcalculi on the
``homomorphic'' and ``antiholomorphic'' subalgebras of $\Sq^{2N-1}$, i.\,e.\
the subalgebras
generated by ${\{z_i~|~i=1,\dots,N\}}$ and ${\{z^*_i~|~i=1,\dots,N\}}$, resp.
\end{itemize}

To end this section, we just mention three further
differential $*$~calculi from Theorem~\ref{fodcf} for which
the bimodule structure takes a simpler form. These are
\begin{itemize}
\item[$\circ$] the calculus $\Gammatpp[q^{2N+2}]$, displaying the simplest
bimodule structure of all calculi:
\[\begin{array}{r@{~}l@{\qquad}r@{~}l}
\d z_kz_l&=q^{-1}\Rm stkl z_s\d z_t;&
\d z_kz^*_l&=q\Rl stkl z^*_s\d z_t;\\[0.5ex]
\d z^*_kz^*_l&=q\Rc stkl z^*_s\d z^*_t;&
\d z^*_kz_l&=q^{-1}\Rrm stkl z_s\d z^*_t.
\end{array}\]
\item[$\circ$] the calculus $\Gammatp[q^4]\equiv\Gammatpp[q^4]$ with
\[\begin{array}{r@{~}l@{\qquad}r@{~}l}
\d z_kz_l&=q^{-1}\Rm stkl z_s\d z_t;&
\d z_kz^*_l&=q\Rlm stkl z^*_s\d z_t-q^{-1}Qz_kz^*_l\H+;\\[0.5ex]
\d z^*_kz^*_l&=q\Rc stkl z^*_s\d z^*_t;&
\d z^*_kz_l&=q^{-1}\Rr stkl z_s\d z^*_t-q^{-1}Qz^*_kz_l\H+.
\end{array}\]
\item[$\circ$] the calculus $\Gammatp[q^2]$ which is also obtained from
$\Gammatl$ in the limit ${\lambda\to q^2}$. Its bimodule structure is
\[\begin{array}{r@{~}l@{\qquad}r@{~}l}
\d z_kz_l&=q^{-1}\Rm stkl z_s\d z_t-q^{-1}Qz_kz_l\H+;&
\d z_kz^*_l&=q\Rlm stkl z^*_s\d z_t;\\[0.5ex]
\d z^*_kz^*_l&=q\Rc stkl z^*_s\d z^*_t-qQz^*_kz^*_l\H+;&
\d z^*_kz_l&=q^{-1}\Rr stkl z_s\d z^*_t.
\end{array}\]
\end{itemize}

None of these calculi is inner. The first two of them decompose
into subcalculi on the holomorphic and antiholomorphic subalgebras.

Note that $\Gammatpd[q^4]\equiv\Gammatppd[q^4]$ is the only isomorphy
of two calculi from Theorem~\ref{fodcf}.

\subsection{Proof of the classification theorems}

\subsubsection{Ansatz obtained by morphisms of tensor products}
We turn now to prove the Theorems~\ref{fodc} and~\ref{fodcf}.
Our approach is based on
an investigation of intertwining mappings for the corepresentations of
the quantum group $\mathcal{A}=\SUq(N)$ on the quantum spheres $X=\Sq^{2N-1}$.
By a similar approach covariant first order differential calculi
on Podle{\'s}' spheres $\mathrm{S}^2_{qc}$ have been classified
by Apel and Schm{\"u}dgen \cite{ap}. Note that $q$ is not a root of unity;
so the representation theory is similar as in the classical case $q=1$.

Let $V(k)$ be the vector space of all $k$-th order polynomials in
the generators $z_i, z^*_i$ of $X$. By the relations~(\ref{qsr1})
and~(\ref{qsr23}), some of these polynomials are identified with
polynomials of lower order. They form a vector subspace in $V(k)$.
Let $\tilde{V}(k)$ denote the complement of this subspace in $V(k)$.
From the coaction $\DeltaR$ we then obtain corepresentations $\pi(k)$ of the
quantum group ${\cal A}= \SUq(N)$ on $X=\Sq^{2N-1}$, with
$\pi(k):\tilde{V}(k)\to \tilde{V}(k)\otimes{\cal A}$.

Here, $\pi(1)$ is the following sum of two irreducible corepresentations,
namely the fundamental representation $u$ of $\SUq(N)$ and its contragredient
$u^{\mathrm{c}}$.

In order to find covariant first order differential calculi on $X$,
intertwining mappings
${T\in\mathrm{Mor}(\pi(1)\otimes\pi(1),\pi(k)\otimes\pi(1))}$
have to be investigated.

To this goal, the direct sum decompositions of the tensor products
${\pi(k)\otimes\pi(1)}$, $k=1,2,\dots$, are calculated.
This can be accomplished e.\,g.\ using Young tableaux; note that
$\pi(k+1)$ is obtained from ${\pi(k)\otimes\pi(1)}$ by removing certain
direct summands according to the commutation relations~(\ref{qsr1}).

Intertwining mappings
${T\in\mathrm{Mor}(\pi(1)\otimes\pi(1),\pi(k)\otimes\pi(1))}$
must correspond to identical direct summands occurring
in both the decompositions of ${\pi(k)\otimes\pi(1)}$,
and ${\pi(1)\otimes\pi(1)}$.
If $N\ge4$, such common summands exist only for $k=1$ (trivial) and $k=3$;
they lead to morphisms from ${\tilde{V}(1)\otimes \tilde{V}(1)}$
to ${\tilde{V}(1)\otimes \tilde{V}(1)}$
and ${\tilde{V}(3)\otimes \tilde{V}(1)}$ which are listed below.
If $N=2$ or $N=3$, there are additional morphisms for ${k=N-1}$,
so the completeness statements in both theorems are guaranteed only for
${N\ge4}$.

The resulting morphisms from ${\tilde{V}(1)\otimes \tilde{V}(1)}$
to ${X\otimes \tilde{V}(1)}$ are given by
\[
\begin{array}{r@{\,\mapsto\,}l@{\quad}r@{\,\mapsto\,}l}
\hbox to2cm{\hfil$z_k\otimes z_l$} &
 \hbox to2cm{$z_k \otimes z_l$\hfil} &
\hbox to2cm{\hfil$z_k\otimes z_l$} &
 \hbox to2cm{$\Rm stkl z_s \otimes z_t$\hfil}\\
z^*_k\otimes z^*_l &  z^*_k\otimes z^*_l&
z^*_k\otimes z^*_l &\Rc stkl  z^*_s\otimes z^*_t\\
z_k\otimes z^*_l & z_k\otimes z^*_l&
z_k\otimes z^*_l &\Rlm stkl z^*_s\otimes z_t\\
z^*_k\otimes z_l &  z^*_k\otimes z_l&
z^*_k\otimes z_l &\Rr stkl  z_s\otimes z^*_t
\end{array}
\]\[
\begin{array}{r@{\,\mapsto\,}l@{\quad}r@{\,\mapsto\,}l}
\hbox to2cm{\hfil$z_k\otimes z^*_l$} &
 \hbox to2cm{$\delta_{kl}q^{2k}\Psi_+$\hfil}&
\hbox to2cm{\hfil$z_k\otimes z^*_l$} &
 \hbox to2cm{$\delta_{kl}q^{2k}\Psi_-$\hfil}\\
z^*_k\otimes z_l & \delta_{kl}\Psi_+&
z^*_k\otimes z_l & \delta_{kl}\Psi_-
\end{array}
\]\[
\begin{array}{r@{\,\mapsto\,}l@{\quad}r@{\,\mapsto\,}l}
\hbox to2cm{\hfil$z_k\otimes z_l$} &
 \hbox to2cm{$z_kz_l\Psi_+$\hfil} &
\hbox to2cm{\hfil$z_k\otimes z_l$} &
 \hbox to2cm{$z_kz_l\Psi_-$\hfil} \\
z^*_k\otimes z^*_l & z^*_k z^*_l\Psi_+&
z^*_k\otimes z^*_l & z^*_k z^*_l\Psi_-\\
z_k\otimes z^*_l &z_k z^*_l\Psi_+&
z_k\otimes z^*_l &z_k z^*_l\Psi_-\\
z^*_k\otimes z_l & z^*_k z_l\Psi_+&
z^*_k\otimes z_l & z^*_k z_l\Psi_-
\end{array}
\]
where $\Psi_+=\sumR{i}z_i\otimes z^*_i$, $\Psi_-=\sumL{i}z^*_i\otimes z_i$.

This leads us to the following ansatz for the bimodule structure of
any covariant differential calculus on $\Sq^{2N-1}$ which is freely
generated as a left module by ${\{\d z_i,\d z^*_i~|~i=1,\dots,N\}}$:
\begin{equation}\label{morans}
\begin{array}{rcl}
\d z_kz_l&=&a_1\Rm stkl z_s\d z_t+b_1z_k\d z_l
  +c_1z_kz_l\H++e_1z_kz_l\H-\\[.7ex]
\d z^*_kz^*_l&=&a_2\Rc stkl z^*_s\d z^*_t+b_2z^*_k\d z^*_l
  +c_2z^*_kz^*_l\H++e_2  z^*_k z^*_l\H-\\[.7ex]
\d z_kz^*_l&=&a_3\Rlm stkl z^*_s\d z_t+b_3z_k\d z^*_l\\*&&
  +c_3z_kz^*_l\H++d_3q^{2k}\delta_{kl}\H+
  +e_3z_kz^*_l\H-+f_3q^{2k}\delta_{kl}\H-\\[.7ex]
\d z^*_kz_l&=&a_4\Rr stkl z_s\d z^*_t+b_4z^*_k\d z_l\\*&&
  +c_4z^*_kz_l\H++d_4\delta_{kl}\H++e_4z^*_kz_l\H-+f_4\delta_{kl}\H-
\end{array}
\end{equation}
with $\H+$ and $\H-$ as defined in section~\ref{srs1}. Here, the $20$ variables
$a_\nu$,~$b_\nu$, $c_\nu$, $d_{\nu^{\prime}}$, $e_\nu$,~$f_{\nu^{\prime}}$
($\nu=1,2,3,4$; $\nu^{\prime}=3,4$)
denote unknown (complex) coefficients.

For the case of covariant first order differential calculi with one
relation ${\H0=0}$, $\H0$ invariant, the ansatz simplifies to
\begin{equation}\label{moransf}
\begin{array}{rcl}
\d z_kz_l&=&a_1\Rm stkl z_s\d z_t+b_1z_k\d z_l
  +c_1z_kz_l\H1\\[.7ex]
\d z^*_kz^*_l&=&a_2\Rc stkl z^*_s\d z^*_t+b_2z^*_k\d z^*_l
  +c_2z^*_kz^*_l\H1\\[.7ex]
\d z_kz^*_l&=&a_3\Rlm stkl z^*_s\d z_t+b_3z_k\d z^*_l
  +c_3z_kz^*_l\H1+d_3q^{2k}\delta_{kl}\H1\\[.7ex]
\d z^*_kz_l&=&a_4\Rr stkl z_s\d z^*_t+b_4z^*_k\d z_l
  +c_4z^*_kz_l\H1+d_4\delta_{kl}\H1
\end{array}
\end{equation}
where $\H1$ is an invariant one-form linearly independent on $\H0$.

\subsubsection{Conditions for the coefficients of the ansatz}
The defining relations of $\Sq^{2N-1}$ together with the properties
required for a covariant first order differential $*$~calculus
can now be used to compile a system of necessary conditions for the
mapping $\d:X\to\Gamma$, and thereby for the coefficients of~(\ref{morans})
and~(\ref{moransf}).

Coefficient comparison in the left module $\Gamma$ is essential for the
following arguments. This yields no difficulty as long as we deal with
differential calculi that are freely generated as left modules, but
needs additional justification in case of the setting of Theorem~\ref{fodcf}.
Now our second classification constraint requires any relation between
one-forms
to be generated algebraically by ${\H0=0}$, where $\H0$ is a linear
combination of $\H+$ and $\H-$; thus, any relation between one-forms
possibly obstructing coefficient comparison needs to involve at least
$N$ of the generators. In fact, none of the coefficient comparisons
done in the following involves more than $3$ independent generators
$\d z_i$ or $\d z^*_i$, so the validity of the classification is
guaranteed for ${N\ge4}$.

\begin{enumerate}
\item\label{fc1} If $\sum\limits_i z_{k_i}^{(*)}z_{l_i}^{(*)}=0$
is one of the (homogeneous) defining relations~(\ref{qsr1}),
$\d\left(\sum\limits_i z_{k_i}^{(*)}z_{l_i}^{(*)}\right)$ must vanish.
Using Leibniz' rule and~(\ref{morans}), this expression can be
written as an element of the left $X$ module generated by~$\d z_i$
and $\d z^*_i$. Since ${\{\d z_i,\d z^*_i\}}$ is required to be a left
module basis, coefficient comparison can be applied.
\item\label{fc2} If $\sum\limits_i z_{k_i}^{(*)}z_{l_i}^{(*)}=0$
is, again, one of the relations~(\ref{qsr1}) and $\d z^{(*)}_m$ any
of the $2N$ bimodule generators, the expression
$\d z_m^{(*)}\,\sum\limits_i z_{k_i}^{(*)}z_{l_i}^{(*)}$ must vanish.
Use again~(\ref{morans}) to write this as a left module expression and
compare coefficients.
\item\label{fc3} The same procedure can be applied to
\begin{equation}\label{dcc3}
\d z_m^{(*)}-\d z_m^{(*)}\left(\sum\limits_i z_i z^*_i\right)
\quad\textrm{and}\quad
q^{-2}\d z_m^{(*)}-\d z_m^{(*)}\left(\sum\limits_i q^{-2i} z^*_i z_i\right)
\end{equation}
both of which must be zero because of~(\ref{qsr23}).
\item\label{fc4} From the $*$~calculus requirement one can
infer $\H++\H+^*=0$. From the definition which gives explicitly a
left module expression for $\H+$ one obtains by the $*$~requirement
a right module expression for $\H+^*$. Rewrite $\H+^*$ to a left module
expression using~(\ref{morans}) and compare coefficients.

Even more conditions can be derived from the $*$~calculus requirement
by taking some expression like $z^{(*)}_k\d z^{(*)}_l$, apply $*$,
(\ref{morans}), $*$, and (\ref{morans}) again, and compare coefficients
with identity.
\end{enumerate}

\subsubsection{The case of freely generated $*$~calculi}
Exploiting conditions from the list given in the preceding section,
we obtain a system of equations for
$a_\nu$,~$b_\nu$, $c_\nu$, $d_{\nu^{\prime}}$, $e_\nu$,~$f_{\nu^{\prime}}$
of the ansatz~(\ref{morans}) which is listed below.
Note that the coefficient comparisons have been done completely for
\ref{fc1} and \ref{fc2} but only in part for \ref{fc3} and \ref{fc4}.
Complex conjugates occur in equations~(\ref{s1})--(\ref{s3}) which
are results of $*$~calculus conditions.

\begin{equation}
\begin{array}{@{}r@{~}l@{}}
a_1&=q^{-1}(b_1+1);\\a_2&=q(b_2+1);\\a_3&=q(b_4+1);\\a_4&=q^{-1}(b_3+1);
\end{array}\quad
\begin{array}{@{}r@{~}l@{}}
c_3&=c_4;\\d_3&=q^{-2N}d_4;\\e_3&=e_4;\\f_3&=q^{-2N}f_4;
\end{array}\quad
\begin{array}{@{}r@{~}l@{}}
a_1a_3&=1\\ a_2a_4&=1\\
qa_4+c_3+q^2\s+ d_3&=0\\
qa_3+e_3+q^2\s+ f_3&=0
\end{array}
\label{al}
\end{equation}
\begin{eqnarray}
b_4(b_3+c_3)+q^{-2}b_1e_3&=&0 \label{b1}\\
(q^{-1}a_1-qa_4+b_1+q^{-2}e_1)b_3+b_2c_1&=&0 \label{b2}\\
q^{-2}e_2b_3+b_2(Qa_2+b_2+c_2)&=&0 \label{b3}\\
b_4c_1+b_1(-Qa_1+b_1+q^{-2}e_1)&=&0 \label{b4}\\
b_3(b_4+q^{-2}e_4)+b_2c_4&=&0 \label{b5}\\
(qa_2+b_2+c_2-q^{-1}a_3)b_4+q^{-2}b_1e_2&=&0 \label{b6}\\
(-q^{-1}a_2+b_2+q^{-1}a_3+q^{-2}e_3)b_3+b_2c_3&=&0 \label{b7}\\
(-qa_1+b_1+qa_4+c_4)b_4+q^{-2}b_1e_4&=&0 \label{b8}\\[.5ex]
(q^{-1}a_3+q^{-2}e_3-qa_2-c_2)c_1 -q^{-2}(c_3+q^2d_3)e_1 \qquad \nonumber \\*
   +b_3c_4 +(qa_4+c_4+d_4-q^{-1}a_1-b_1)c_3 &=&0 \label{c1}\\
(b_3+e_3)e_4+c_3f_4-b_1e_3-c_1e_2-e_1f_3&=&0 \label{c2}\\
b_4c_3-b_2c_4+q^{-2}(c_3+q^2d_3)e_4-d_4c_2-q^{-2}c_1e_2&=&0 \label{c3}\\
(qa_4+c_4-q^{-1}a_1-q^{-2}e_1)e_2 -c_2(e_4+f_4) \qquad \nonumber \\*
   +(q^{-1}a_3+q^{-2}e_3+f_3-q^{-1}a_2-b_2)e_4&=&0 \label{c4}\\[.5ex]
(-q^{-1}a_1+qa_4+c_4+d_4)d_3+q^{-2}c_1f_3&=&0 \label{d1}\\
(e_4+f_4)d_3+q^{-2}e_1f_3&=&0 \label{d2}\\
c_2d_4+q^{-2}(c_3+d_3)f_4&=&0 \label{d3}\\
e_2d_4+q^{-2}(-q^3a_2+qa_3+e_3+f_3)f_4&=&0 \label{d4}\\
(Qa_4+c_4+d_4)d_4+q^{-2}c_1f_4&=&0 \label{d5}\\
(e_4+f_4)d_4+q^{-2}(qa_1+e_1-qa_4)f_4&=&0 \label{d6}\\
(qa_2-qa_3+c_2)d_3+q^{-2}(c_3+q^2d_3)f_3&=&0 \label{d7}\\
e_2d_3+q^{-2}(-q^2Qa_3+e_3+q^2f_3)f_3&=&0 \label{d8}\\[.5ex]
b_4d_3+q^{-2}b_1f_3&=&0 \label{e1}\\
b_2d_4+q^{-2}b_3f_4&=&0 \label{e2}\\[.5ex]
\cc{a_3}b_4+a_3(\cc{b_3}-\cc{c_3}b_4-q^{-2}\cc{e_3}b_1)&=&0 \label{s1}\\
\cc{a_4}b_3+a_4 (\cc{b_4}-\cc{c_4}b_2-q^{-2}\cc{e_4}b_3)&=&0 \label{s4}\\
q\cc{a_4}e_3+e_4((\cc{b_4}-\cc{c_4}b_2-q^{-2}\cc{e_4}b_3)
   \qquad\qquad\qquad\qquad \nonumber \\
   -(qa_4+c_4+d_4)(\cc{c_4}(qa_2+c_2)+q^{-2}\cc{e_4}(c_3+q^2d_3))
   \qquad \nonumber \\
   -q^{-2}c_1(\cc{c_4}e_2+q^{-2}\cc{e_4}(qa_3+e_3+q^2f_3))&=&0 \label{s5}\\
(q^{2N+1}\cc{a_4}+\cc{b_4}-\cc{c_4}b_2-q^{-2}\cc{e_4}b_3)d_3-\cc{d_3}&=&0
   \label{s2}\\
(q^{2N+1}\cc{a_4}+\cc{b_4}-\cc{c_4}b_2-q^{-2}\cc{e_4}b_3)f_3-\cc{f_3}&=&0
   \label{s3}
\end{eqnarray}

We stress that this (partially redundant) system of equations
gives a set of necessary conditions which is not a priori
complete since not all required properties of a first order
covariant differential $*$~calculus have been fully exploited.
So, for any set of coefficients solving this system,
it remains still necessary to prove that all required properties
are satisfied.

In order to solve this system of equations, we observe first
that~(\ref{d1})--(\ref{d3}) together with~(\ref{al}) imply that
\begin{equation}
d_3=d_4=0~~\mbox{if and only if}~~f_3=f_4=0.
\label{eqvd}
\end{equation}
Furthermore, equation~(\ref{s1}) together with~(\ref{al}) implies that
\begin{equation}
b_1=0~~\mbox{if and only if}~~b_2=0.
\label{eqvb}
\end{equation}
Note that~(\ref{eqvb}), unlike~(\ref{eqvd}), depends on the $*$~calculus
requirement.

Assume now that ${b_1\ne0}$ and ${d_3\ne0}$.
In this case, we can express $a_{\nu}$, $b_{\nu}$, $c_{\nu'}$, $d_{\nu'}$,
$e_{\nu'}$, $f_{\nu'}$ (${\nu=1,2,3,4}$; ${\nu'=3,4}$) in terms of
only two parameters ${\alpha:=q^{-1}a_1}$ and ${\tau:=f_3}$, using
the equations~(\ref{al}), (\ref{e1}),
and~(\ref{e2}). Subsequently, $c_1$, $c_2$, $e_1$, and $e_2$ are also
expressed in terms of these parameters by using~(\ref{d1})--(\ref{d4}).

Let now ${d_3=0}$. Then the equations~(\ref{e1}) and~(\ref{e2}) don't lead
to any additional restriction. Instead, we use~(\ref{al}),
(\ref{d1})--(\ref{d4}), (\ref{s4}), and (\ref{s5}) to express
$a_{\nu}$, $b_{\nu}$, $c_{\nu}$, $e_{\nu}$ by ${\alpha=q^{-1}a_1}$,
${\omega=c_1}$ and ${\psi=c_2+\calpha^{-1}}$. Then, for ${b_1\ne0}$
equation~(\ref{b3}) implies ${\psi=q^2\calpha\omega^{-1}+1}$, yielding
the bimodule structure of~$\Gammaapd$. For ${b_1=0}$, $\Gammapd$ is
obtained.

Finally, if ${b_1=0}$ and ${d_3\ne0}$, the coefficients $a_{\nu}$ and
$b_{\nu}$ are already determined while $c_{\nu}$, $e_{\nu}$, $d_{\nu'}$,
$f_{\nu'}$ (${\nu=1,2,3,4}$; ${\nu'=3,4}$) depend on two parameters
${\varrho=-q^2d_3}$, ${\tau=-f_3}$ by equations~(\ref{al}),
(\ref{d1})--(\ref{d4}), yielding the bimodule structure of~$\Gammappd$.

For all calculi, $*$~conditions imply that the parameters be real, except
for $\Gammapd$ where $\alpha$ can take non-real values if the additional
condition ${\omega=q^4\alpha\calpha}$ is fulfilled.

One checks that $\Gammaatd$, $\Gammaapd$, $\Gammapd$, $\Gammappd$
satisfy all requirements for a covariant first order differential
$*$~calculus:
Covariance is guaranteed by the ansatz~(\ref{morans}).
The first two groups of conditions in the list above resulting from the
relations~(\ref{qsr1}) are completely encoded in the system of equations,
so these are also fulfilled.
Finally, it is checked by direct calculations that the elements~(\ref{dcc3})
vanish and the $*$~calculus property is satisfied.
This completes the proof of Theorem~\ref{fodc}.\qed

\subsubsection{Calculi with the relation $\H0=0$}

We shall now use again the conditions of types~\ref{fc1}--\ref{fc3}
in order to specify the possible sets of coefficients in~(\ref{moransf})
for first order differential calculi with one relation
\begin{equation}
\H++\lambda\H-=0,\quad\lambda\ne0.
\label{Hl}
\end{equation}

Much as in the case of freely generated calculi, we evaluate
the conditions~\ref{fc1}, observing now the additional condition~(\ref{Hl}),
to obtain the equations
\begin{eqnarray}
&\begin{array}{@{}r@{~}l@{}}
a_1&=q^{-1}(b_1+1);\\a_2&=q(b_2+1);
\end{array}\quad
\begin{array}{@{}r@{~}l@{}}
a_3&=q(b_4+1);\\
a_4&=q^{-1}(b_3+1);
\end{array}\quad
\begin{array}{@{}r@{~}l@{}}
c_3&=c_4;\\d_3&=q^{-2N}d_4;
\end{array}\quad
\begin{array}{@{}r@{~}l@{}}
a_1a_3&=1;\\ a_2a_4&=1;
\end{array}&
\label{alf}
\\
&-q\lambda^{-1}a_3+b_3+c_3+q^2\s+d_3=-1.&
\label{alfi}
\end{eqnarray}
The condition~(\ref{Hl}) implies ${(\H++\lambda\H-)z^{(*)}_m=0}$ which
leads to
\begin{eqnarray}
q^{-2}\lambda b_1+b_4&=&0 \label{f1}\\*
q^2\lambda^{-1}b_2+b_3&=&0 \label{f2}\\*
(qa_4+c_4+d_4-q^{-1}a_1)+q^{-2}\lambda c_1&=&0 \label{f3}\\
(qa_2+c_2-q^{-1}a_3)+q^{-2}\lambda(c_3+q^2d_3)&=&0. \label{f4}
\end{eqnarray}
From~(\ref{f1}), (\ref{f2}) it follows by means of~(\ref{alf}) that
\begin{eqnarray}
&(a_1=q\lambda^{-1}~~\mbox{or}~~a_1=q^{-1})& \label{f5} \\
&(a_2=q^{-1}\lambda~~\mbox{or}~~a_2=q).& \label{f6}
\end{eqnarray}
By virtue of~(\ref{f3}), (\ref{f4}), (\ref{alf}) and~(\ref{alfi}),
$c_\nu$ can be expressed in terms of $a_1$, $a_2$, $d_3$.
Evaluation of conditions
of type~\ref{fc2} then yields that $d_3$ has to be zero, except if
${a_1=q^{-1}}$ and ${a_2=q}$. In the latter case, $d_3$ can still take
the values $0$ or ${q^{-2}\se+^{-1}(q^4\lambda^{-1}-1)}$.
By combination, we obtain five cases which give the differential calculi
$\Gammatld$, $\Gammatpd$, $\Gammatppd$ (${\lambda\ne0,\infty}$),
$\Gammatod$, $\Gammatood$ from Theorem~\ref{fodcf}, resp.

Similar considerations based on the relations ${\H+=0}$ or ${\H-=0}$
instead of~(\ref{Hl}) lead to $\Gammatppd[\infty]$ and $\Gammatppd[0]$,
resp.

Again, one checks all required properties to establish that all of these
are first order differential calculi, and that only the three series
$\Gammatld$, $\Gammatpd$, $\Gammatppd$ for real $\lambda$ are $*$~calculi.
Thus, Theorem~\ref{fodcf} is proved.\qed

\section{Higher order calculus and symmetry}\label{hos}
\subsection{Higher order differential calculi}
In order to describe the noncommutative differential geometry of
quantum spaces, it is insufficient to have only first order
differential calculus since basic concepts of differential geometry
require at least second order differential forms. This causes us
to turn our attention to higher order differential calculi on
quantum homogeneous spaces.

Our definitions for higher order differential calculi on quantum
homogeneous spaces follow those given in~\cite{wr} for the bicovariant case
in the quantum group setting.

Let $X$ a quantum homogeneous space for the quantum group $\mathcal{A}$.
Let $(\Gamma,\d)$ be a covariant first order differential
calculus on~$X$.

Then, a \emph{covariant higher order differential calculus} on $X$ is a pair
$(\Gamma^{\wedge},\d)$ consisting of a graded algebra $\Gamma^{\wedge}$ with
multiplication denoted by~$\wedge$ and a linear mapping
$\d:\Gamma^{\wedge}\to\Gamma^{\wedge}$ such that the following conditions
are fulfilled:
\begin{enumerate}
\item The degree~$0$ and degree~$1$ components of~$\Gamma^{\wedge}$
are isomorphic to~$X$ and~$\Gamma$, respectively (they will be identified
with $X$ and $\Gamma$ in the following).
\item The mapping~$\d$ increases the degree by~$1$, and
$\d$ extends the differential~$\d:X\to\Gamma$ from the first
order differential calculus $(\Gamma,\d)$.
\item The mapping~$\d$ is a \emph{graded derivative,} i.\,e.\ it fulfils the
graded Leibniz' rule
\[\d(\vartheta_1\wedge\vartheta_2)=\d\vartheta_1\wedge\vartheta_2
+(-1)^d\vartheta_1\wedge\d\vartheta_2\]
for any ${\vartheta_1,\vartheta_2\in\Gamma^{\wedge}}$
(with $d$ being the degree of $\vartheta_1$).
For any $\vartheta\in\Gamma^{\wedge}$, $\d(\d\vartheta)=0$.
\item The covariance map ${\PhiR:\Gamma\to\Gamma\otimes{\cal A}}$
from the first order differential calculus can be extended to a map
$\PhiR^{\wedge}:\Gamma^{\wedge}\to\Gamma^{\wedge}\otimes{\cal A}$
making $\Gamma^{\wedge}$ into a covariant $X$-bimodule.
\end{enumerate}

If the underlying first order calculus ${(\Gamma,\d)}$ is a $*$~calculus,
then there is an induced $*$~structure on ${(\Gamma^{\wedge},\d)}$,
and $\d(\vartheta^*)=(\d\vartheta)^*$ is fulfilled
for all~$\vartheta\in\Gamma^{\wedge}$.

\subsection{Higher order differential calculi on quantum spheres}
Now we turn again to study the quantum spheres $\Sq^{2N-1}$.
Since our main interest is to provide a framework of covariant differential
calculus appropriate to describe the noncommutative geometry of
the quantum spheres, our considerations should be based on one of
the differential calculi from Theorem~\ref{fodcf} that are not freely
generated (as left modules) by the generator set
${\{\d z_i,\d z^*_i~|~i=1,\dots,N\}}$ but carry an additional relation
of the type ${\H++\lambda\H-=0}$, as holds in the classical case ${q=1}$ with
${\lambda=1}$.

We restrict our considerations from now on to
the differential calculus ${\Gammad\equiv\Gammatld[1]}$ with the
relation~\ref{fr} and bimodule structure~\ref{fodcfb} given in~\ref{srs1}.
Since we want to transfer ideas and techniques from the theory of
bicovariant differential calculi on quantum groups, it is essential
to work with an inner calculus since all bicovariant calculi are inner
in the quantum group case;
on the other hand, for many of the calculations done
in this section we need the decomposition of $\Gammat$ into subcalculi
on the holomorphic and antiholomorphic subalgebras of ${X=\Sq^{2N-1}}$.

As a consequence of~(\ref{fr}), we have the relation $\d\H++\d\H-=0$ or
\begin{equation}\label{dh}
\sumR{i}\d z_i\wedge\d z^*_i+\sumL{i}\d z^*_i\wedge\d z_i=0
\end{equation}
in any higher order differential calculus extending ${(\Gamma,\d)}$.
The bimodule structure of ${(\Gamma,\d)}$ implies
\begin{equation}\label{hh}
\d\H{}=\H{}\wedge\H{}
\end{equation}
(remember that $\H{}=qQ^{-1}\H-$, thus $\d\H{}=qQ^{-1}\d\H-$),
and the following commutation relations for
the generators $\d z_i$, $\d z^*_i$:
\begin{equation}\label{hodcr}
\begin{array}{rcl}
0 &=& \d z_k\wedge\d z_l+q\R stkl \d z_s\wedge\d z_t\\[.7ex]
0 &=& \d z^*_k\wedge\d z^*_l+q^{-1}\Rcm stkl \d z^*_s\wedge\d z^*_t\\[.7ex]
0 &=&
q^{-3}\Rlm stkl(\d z^*_s\wedge\d z_t+Q^2z^*_s\d z_t\wedge\H{})
\\*&&\quad
+(\d z_k\wedge\d z^*_l+Q^2z_k\d z^*_l\wedge\H{})
+Q^2(Q^2+1)z_k z^*_l\d\H{}.
\end{array}
\end{equation}

There is a universal higher order differential calculus
${(\Gamma^{\wedge}_{\mathrm{u}},\d)}$
which is generated as an algebra by the $4N$ elements
$z_i$, $z^*_i$, $\d z_i$, $\d z^*_i$
subject to the relations~(\ref{qsr1}), (\ref{qsr23}), (\ref{fr}),
(\ref{fodcfb}), (\ref{dh}),~(\ref{hh}),~(\ref{hodcr}).

The calculus ${(\Gamma^{\wedge}_{\mathrm{u}},\d)}$ is not an
\lq\lq inner\rq\rq\ calculus in the sense of the differential mapping
$\d$ being generated by a graded commutator.
So by imposing this as an additional condition, i.\,e.\
\begin{equation}\label{gcr}
\d\vartheta=\H{}\wedge\vartheta-(-1)^d\vartheta\wedge\H{},\quad
d=\textrm{degree of $\vartheta$},
\end{equation}
for all $\vartheta$,
we obtain a smaller higher order differential calculus which
we will denote by ${(\Gamma^{\wedge}_*,\d)}$.
Obviously, (\ref{dh}) and (\ref{gcr}) together imply
\[\d\H{}=0.\]

\begin{Prop}
In ${(\Gamma^{\wedge}_*,\d)}$, all differential forms of degree
$2N-1$ or higher vanish.

If $N\ge3$, there is one differential $(2N-2)$-form which generates
the set of all $(2N-2)$-forms as a left module.
\end{Prop}

For the proof we need to consider differential forms
\[\Theta=\d z_{j_1}\wedge\dots\wedge\d z_{j_\mu}\wedge
\d z^*_{k_1}\wedge \dots\wedge\d z^*_{k_\nu}\]
with $1\le j_1<\dots<j_\mu\le N$, $1\le k_1<\dots<k_\nu\le N$.
These elements generate $\Gamma^{\wedge}_*$ as a left module.

Then, the main argument used to prove both parts of the Proposition is
stated in the following lemma:
\begin{Lemma}
If $\Theta$ is chosen as above and if the index sets
$J=\{j_1,\dots,j_\mu\}$ and $K=\{k_1,\dots,k_\nu\}$ fulfil
$J\cup K=\{1,\dots,N\}$ and $J\cap K\ne\emptyset$, then $\Theta=0$.
\end{Lemma}

\begin{Proof}\textbf{of the Lemma~}
Note that all $\d z_i$ quasi-commute, as do all $\d z^*_i$.
Therefore we can use the substitution
\[\d z_j\wedge\d z^*_j=-\sum\limits_{i\ne j}\d z_j\wedge\d z^*_j\]
(resulting from $\d\H{}=0$) to rewrite $\Theta$ as a sum of $N-1$ members
each of which contains one of the expressions $\d z_i\wedge\d z_i$,
$\d z^*_i\wedge\d z^*_i$ with $1\le i\le N$.
\end{Proof}

\begin{Proof}\textbf{of the Proposition~}
It is easily seen that each differential form of degree $2N-1$
is a linear combination (with coefficients from $\Sq^{2N-1}$)
of forms of type $\Theta$ fulfilling the additional index set condition,
and thus vanishes. This proves the first part of the Proposition.

Moreover, for any differential $(2N-2)$-form $\Theta$ of the above type
it can be seen that the index sets $J$ and $K$ either fulfil the
index set condition leading to $\Theta=0$, or
$J=K=\{1,\dots,N\}\setminus\{j\}$ with a single index $j\in\{1,\dots,N\}$.
But all $\Theta$'s with the latter property are transformed into
(scalar) multiples of each other by just the same substitution
as in the proof of the lemma. This completes the
proof of the second part.
\end{Proof}

\subsection{A symmetry concept for quantum spaces}\label{as}
For quantum groups, Woronowicz~\cite{wr} has described a construction
extending a bicovariant first order differential calculus~$(\Gamma, \d)$
to a bicovariant higher order differential calculus~${(\Gamma^{\wedge},\d)}$
by antisymmetrisation. For the antisymmetrisation procedure,
a bimodule homomorphism of the tensor product $\Gamma\otimes\Gamma$
is required,%
\footnote{Here and in the following tensor products of differential
modules are always meant to be tensor products over the corresponding
quantum group or space, i.\,e.\ ${\Gamma\otimes\Gamma}$ means
${\Gamma\otimes_{\mathcal{A}}\Gamma}$ or ${\Gamma\otimes_X\Gamma}$,
and so on}
$\sigma:\Gamma\otimes\Gamma\to\Gamma\otimes\Gamma$,
which needs to fulfil the following \emph{braid equation} on the three-fold
tensor product $\Gamma\otimes\Gamma\otimes\Gamma$:
\begin{equation}
\label{sb}
(\sigma\otimes\id)\circ(\id\otimes\sigma)\circ(\sigma\otimes\id)\equiv
(\id\otimes\sigma)\circ(\sigma\otimes\id)\circ(\id\otimes\sigma)
\end{equation}
or, in short,
$\sigma_{12}\sigma_{23}\sigma_{12}=\sigma_{23}\sigma_{12}\sigma_{23}$.
Here, $\sigma_{i\,i+1}$ denote actions of $\sigma$ on the
$i$-th and $(i+1)$-th components of a multiple tensor product.

The antisymmetrisation procedure then works as follows:
For any permutation $p$ of $\{1,\dots,k\}$ let $\ell(p)$ denote the
length of $p$, i.\,e.\ the number of inversions in $p$.
Then $p$ has a decomposition
$p=t_{i_1}\circ t_{i_2}\circ\dots\circ t_{i_\ell}$
where $1\le i_1,\dots,i_{\ell}\le p-1$, and $t_i$ means
the permutation of length~$1$ (transposition) that exchanges the
$i$-th and $(i+1)$-th elements. Let
$\sigma^p=\sigma_{i_1\,i_1+1}\circ\sigma_{i_2\,i_2+1}\circ\dots\circ
\sigma_{i_\ell\,i_\ell+1}$. Because of the braid relation~(\ref{sb}),
$\sigma^p$ is independent on the choice of the decomposition of $p$
and thereby well-defined.
By $A_k=\sum\limits_p(-1)^{\ell(p)}\sigma^p$ an antisymmetriser
on the $k$-fold tensor product $\Gamma^{\otimes k}$
is defined. The $k$-th degree component of the higher order differential
algebra $\Gamma^{\wedge k}$, and the higher order differential algebra
$\Gamma^{\wedge}$ are obtained by
\[\Gamma^{\wedge k}=\Gamma^{\otimes k}/\ker A_k;\qquad
\Gamma^{\wedge}=\sum\limits_{k=0}^\infty\!\!{}^\oplus\,\Gamma^{\wedge k}.\]

Since $\sigma$ encodes a symmetry in $\Gamma^{\otimes}$ which is
also relevant for the noncommutative geometry of the quantum group,
we want to transfer Woronowicz's construction
to our setting of a quantum homogeneous space with (only one-sided)
covariant first order differential calculus.

We start with a definition for a symmetry homomorphism which is already
adapted for a slightly generalised situation, admitting factorisation
of the underlying tensor product---we shall need this later.

\begin{Def}
Let ${\cal A}$ be a quantum group and $X$ a quantum space for
${\cal A}$. Let $(\Gamma,\d)$ be a first order differential calculus
for $X$ which is covariant w.\,r.\,t.~${\cal A}$, and $M$ a sub-bimodule
in $\Gamma\otimes\Gamma$.
Let $\Gamma^{\otimes2}_M:=(\Gamma\otimes\Gamma)/M$.

Then, a \emph{symmetry homomorphism for $\Gamma^{\otimes2}_M$} is a
bimodule homomorphism $\sigma:\Gamma^{\otimes2}_M\to\Gamma^{\otimes2}_M$
for which the braid equation~(\ref{sb}) is fulfilled in
${\Gamma^{\otimes3}_M:=(\Gamma\otimes\Gamma\otimes\Gamma)/M_3}$
where $M_3$ is the closure under $\sigma_{12}$ and $\sigma_{23}$
of the subbimodule generated by
${M\otimes X+X\otimes M}$ in ${\Gamma\otimes\Gamma\otimes\Gamma}$.
\end{Def}
\subsection{Symmetry on quantum spheres}
We are interested in symmetry homomorphisms which lead to non-trivial
higher order differential calculi.
By a non-trivial higher order differential calculus we shall mean, in the
following statements, a higher order differential calculus in which
the differential $2$-forms $\d z_k\wedge\d z_l$ and $\d z^*_k\wedge\d z^*_l$
for $k\ne l$ are nonzero.

First we let $M=\{0\}$, thereby seeking symmetry homomorphisms
on~$\Gamma\otimes\Gamma$.

\begin{Prop}\label{nosy}
On~$\Gamma\otimes\Gamma$, there
is no bimodule homomorphism~$\sigma$ leading to non-trivial higher order
differential calculus.
\end{Prop}

\begin{Proof}
Assume there is a bimodule homomorphism~$\sigma$
on~$\Gamma\otimes\Gamma$
leading to non-trivial higher order differential calculus.
The latter requires $\sigma$ to fulfil (among others) the following
equations:
\begin{eqnarray}
(\sigma-\id)(\d z_k\otimes\d z_l+q\R stkl \d z_s\otimes\d z_t)
  &=& 0   \label{eq-smidoo}\\
(\sigma-\id)(\d z^*_k\otimes\d z^*_l+q^{-1}\Rcm stkl \d z^*_s\otimes\d z^*_t)
  &=& 0.  \label{eq-smidss}
\end{eqnarray}
Equation~(\ref{eq-smidoo}) implies that $\sigma$ is given on elements
of the form $\d z_k\otimes\d z_l$ by one of the following three
equations:
\[\begin{array}{r@{\qquad}rcl@{\qquad}l}
&\sigma(\d z_k\otimes\d z_l)&=&\d z_k\otimes\d z_l&\textbf{(A1)}\\
\textrm{or}&\sigma(\d z_k\otimes\d z_l)&=&q^{-1}\R stkl \d z_s\otimes\d z_t
  &\textbf{(A2)}\\
\textrm{or}&\sigma(\d z_k\otimes\d z_l)&=&q\Rm stkl \d z_s\otimes\d z_t;
  &\textbf{(A3)}
\end{array}\]
equally it follows from equation~(\ref{eq-smidss}) that, on elements
of the form $\d z^*_k\otimes\d z^*_l$, the homomorphism~$\sigma$
is given by one of the three equations:
\[\begin{array}{r@{\qquad}rcl@{\qquad}l}
&\sigma(\d z^*_k\otimes\d z^*_l)&=&\d z^*_k\otimes\d z^*_l&\textbf{(B1)}\\
\textrm{or}
&\sigma(\d z^*_k\otimes\d z^*_l)&=&q\Rcm stkl \d z^*_s\otimes\d z^*_t
  &\textbf{(B2)}\\
\textrm{or}
&\sigma(\d z^*_k\otimes\d z^*_l)&=&q^{-1}\Rc stkl \d z^*_s\otimes\d z^*_t.
  &\textbf{(B3)}
\end{array}\]
Now since (A1) or (B1) would annihilate any $2$-form $\d z_k\wedge\d z_l$,
$\d z^*_k\wedge\d z^*_l$, resp., in the higher order differential calculus,
only the cases (A2) and (A3), (B2) and (B3), resp., need to be considered.

By the identities
\begin{eqnarray*}
\H{}\otimes\d z^*_k &=& -qQ^{-1}\sumR{i} z_i\d z^*_i\otimes\d z^*_k,\\
\H{}\otimes\d z_k   &=& qQ^{-1}\sumL{i} z^*_i\d z_i\otimes\d z_k,
\end{eqnarray*}
and the required left module homomorphism property of $\sigma$,
expressions for $\sigma(\H{}\otimes\d z^{(*)}_k)$ are obtained which
depend on the respective case conditions (A2) or (A3), (B2) or (B3).
From those we infer, by a similar decomposition of $\H{}$ in the second
tensor factor, and the right module homomorphism requirement of $\sigma$,
the following expressions for $\sigma(\HtH)$:
\[\begin{array}{rcll}
\sigma(\HtH)&=&\HtH+
qQ^{-1}\left(\sumR{i}\d z_i\otimes\d z^*_i
-q^{-2}\sumL{i}\d z^*_i\otimes\d z_i\right)\\[.5ex]
&&&\makebox[0pt][r]{(cases (A2), (B3))}\\[1ex]
\sigma(\HtH)&=&\HtH+
q^3Q^{-1}\left(\sumR{i}\d z_i\otimes\d z^*_i
-q^{-2}\sumL{i}\d z^*_i\otimes\d z_i\right)\\[.5ex]
&&&\makebox[0pt][r]{(cases (A3), (B2))}.
\end{array}\]
Since both expressions differ, our case distinction reduces from
four to two possible combinations of cases, (A2)/(B3) and (A3)/(B2).

For these cases, expressions for
${\sigma(\d z^*_k\otimes\d z_l)}$ and ${\sigma(\d z_k\otimes\d z^*_l)}$
can be calculated from ${\sigma(\H{}\otimes\d z_k^{(*)})}$
using~(\ref{foi}), (\ref{fodcfb}), and the bimodule homomorphism
requirement for~$\sigma$;
but finally it turns out that the braid relation is violated.
This is demonstrated by considering one particular element of
$\Gamma\otimes\Gamma\otimes\Gamma$; namely, we have
\[
(\sigma_{12}\sigma_{23}\sigma_{12}-\sigma_{23}\sigma_{12}\sigma_{23})
\left(\sumR{i}\H{}\otimes\d z_i\otimes\d z^*_i\right)
\ne 0.\]
This implies that even under the conditions (A2)/(B3) or (A3)/(B2)
no symmetry homomorphism for ${\Gamma\otimes\Gamma}$ is obtained.
\end{Proof}

\begin{Lemma}
In~$\Gamma\otimes\Gamma$, the elements
\begin{equation}\label{hof}
\sumL{i}\d z^*_i\otimes\d z_i -
q^2\sumR{i}\d z_i\otimes\d z^*_i \textrm{\quad and\quad}
\sumL{i}\d z^*_i\otimes\d z_i + Q^2\HtH
\end{equation}
generate a sub-bimodule $M_0$ as a left module.
\end{Lemma}

\begin{Proof}
Both of these elements quasi-commute with all $x\in X$.
\end{Proof}

Define
\[ \Gamma^{\otimes2}_\circ:=(\Gamma\otimes\Gamma)/{M_0}.\]

\begin{Prop}\label{syf}
There are exactly two symmetry homomorphisms on $\Gamma^{\otimes2}_{\circ}$
leading to non-trivial higher order differential calculi.
They are inverse to each other.

One of them is the bimodule homomorphism $\sigma$ defined by
\begin{equation}\label{sig}
\begin{array}{rcl}
\sigma(\d z_k\otimes\d z_l) &=&
q^{-1}\R stkl \d z_s\otimes\d z_t\\[6pt]
\sigma(\d z^*_k\otimes\d z^*_l) &=&
q^{-1}\Rc stkl \d z^*_s\otimes\d z^*_t\\[6pt]
\sigma(\d z_k\otimes\d z^*_l) &=&
q^{-3}\Rlm stkl \left(\d z^*_s\otimes\d z_t
+Q^2z^*_s\d z_t\otimes\H{}\right)\\*
&&\quad+q^{-1}Q\d z_k\otimes\d z^*_l\\*
&&\quad-Q^2z_k\H{}\otimes\d z^*_l -q^{-3}Q^3z_kz^*_l\HtH\\[6pt]
\sigma(\d z^*_k\otimes\d z_l) &=&
q\Rr stkl \left(\d z_s\otimes\d z^*_t
+Q^2z_s\d z^*_t\otimes\H{}\right)\\*
&&\quad-q^{-2}Q^2z^*_k\H{}\otimes\d z_l+qQ^3z^*_kz_l\HtH.
\end{array}
\end{equation}
\end{Prop}

\begin{Proof}
Since the factorisation of the tensor product brings about no change for
the parts of the tensor product which are generated only by $\d z_i$
resp.\ only by $\d z^*_i$,
we can start as in the proof of Prop.~\ref{nosy} and obtain the same
possible cases (A2), (A3) for ${\sigma(\d z_k\otimes\d z_l)}$, and
(B1), (B3) for ${\sigma(\d z^*_k\otimes\d z^*_l)}$ as before.
The subsequent calculations for ${\sigma(\H{}\otimes\d z^{(*)}_k)}$
and ${\sigma(\HtH)}$ remain valid, too; but the resulting expressions
are now simplified since the invariant elements~(\ref{hof}) are
zero in $\Gamma^{\otimes2}_{\circ}$. In particular,
our argument used above to rule out the combinations
(A2)/(B2) and (A3)/(B3) fails since we have now in all cases
\[\sigma(\HtH)=\HtH.\]

To deal, therefore, with (A3)/(B3)
(the argument is quite the same for the other case)
we calculate expressions for
${\sigma(\d z_k\otimes\d z^*_l)}$ and ${\sigma(\d z^*_k\otimes\d z_l)}$
just as done for (A2)/(B3) in the proof of Prop.~\ref{nosy}.
Then we obtain, for any $k\ne l$,
\[\begin{array}{l}
(\sigma-\id)(q\d z_l\otimes\d z^*_k +q^{-2}\d z^*_k\otimes\d z_l
+qQ^2 z_l\d z^*_k\otimes\H{} +q^{-2}Q^2 z^*_k\d z_l\otimes\H{})
\\[1ex] \quad=
-Q\d z_l\otimes\d z^*_k +q^{-1}Q\d z^*_k\otimes\d z_l -Q^4z^*_kz_l\HtH
\\ \quad\qquad
+qQ^2(z_l\H{}\otimes\d z^*_k-z_l\d z^*_k\otimes\H{})
+q^{-2}Q^2(z^*_k\H{}\otimes\d z_l-z^*_k\d z_l\otimes\H{}).
\end{array}\]
Here, the left-hand side should vanish
for any $\sigma$ leading to non-trivial higher order differential
calculus, but the right-hand side does not.

So we are once more left with (A2)/(B3) and (A3)/(B2). By calculating
${\sigma(\d z_k\otimes\d z^*_l)}$ and ${\sigma(\d z^*_k\otimes\d z_l)}$
for (A2)/(B3), we find the last two equations of the Proposition.

Consider now the left module homomorphism
${\sigma_{\mathrm{L}}:\Gamma\otimes\Gamma}\to{\Gamma\otimes\Gamma}$ and
the right module homomorphism
${\sigma_{\mathrm{R}}:\Gamma\otimes\Gamma}\to{\Gamma\otimes\Gamma}$
defined by the equations~(\ref{sig}).
In order to prove that the bimodule homomorphism $\sigma$ is well-defined
one checks that
\begin{enumerate}
\item the subbimodule $M_0$ of ${\Gamma\otimes\Gamma}$
is invariant under both $\sigma_{\mathrm{L}}$ and $\sigma_{\mathrm{R}}$;
\item for any ${a\in\Gamma\otimes\Gamma}$,
$\sigma_{\mathrm{L}}(a)-\sigma_{\mathrm{R}}(a)\in M_0$.
\end{enumerate}
For the second part it is sufficient to consider expressions
of the type ${\d z^{(*)}_j\otimes\d z^{(*)}_k\cdot z^{(*)}_l}$,
e.\,g.\ for ${\d z_j\otimes\d z_k\cdot z_l}$ (the simplest case) we have
\[
\begin{array}{@{}r@{~}l@{~}l@{}}
\sigma_{\mathrm{L}}(\d z_j\otimes\d z_k\cdot z_l)
&=q^2\R uvjs \R stkl z_u\sigma_{\mathrm{L}}(\d z_v\otimes\d z_t)
&=q\R abvt \R uvjs \R stkl z_u\d z_a\otimes\d z_b\\[.7ex]
\sigma_{\mathrm{R}}(\d z_j\otimes\d z_k\cdot z_l)
&=q^{-1}\R stjk \d z_s\otimes\d z_t\cdot z_l
&=q\R uasv \R vbtl \R stjk z_u\d z_a\otimes\d z_b.
\end{array}
\]

Finally, the braid relation~(\ref{sb}) has to be proved for $\sigma$.
Because of (\ref{foi}) and the bimodule homomorphism property of $\sigma$,
it is sufficient to show~(\ref{sb}) on the elements
$\d z^{(*)}_k\otimes\H{}\otimes\d z^{(*)}_l$. To give again the simplest case,
\[\begin{array}{l}
\sigma_{12}\sigma_{23}\sigma_{12}(\d z_k\otimes\H{}\otimes\d z_l)=
\sigma_{23}\sigma_{12}\sigma_{23}(\d z_k\otimes\H{}\otimes\d z_l)\\[1ex]
\qquad\qquad=q^{-4}Q\R stkl \H{}\otimes\d z_s\otimes\d z_t
+q^{-3}\R stkl \d z_s\otimes\H{}\otimes\d z_t
\end{array}\]

Another symmetry homomorphism $\sigma^{\prime}$ is obtained by applying
the same procedure to the case (A3)/(B2). From the defining equations
for $\sigma$ and $\sigma^{\prime}$ it can be seen that
${\sigma\circ\sigma^{\prime}}\equiv{\sigma^{\prime}\circ\sigma\equiv\id}$.
Obviously, $\sigma\ne\sigma^{\prime}$.
\end{Proof}

As a consequence of Proposition~\ref{syf}, we are now able to transfer
Woronowicz's antisymmetrisation construction to the quantum spheres.
Let ${\Gamma^{\otimes}_{\circ}:=\Gamma^{\otimes}/I_0}$ where
$I_0$ is the closure under actions of $\sigma_{i,i+1}$ of
the two-sided ideal in $\Gamma^{\otimes}$ generated by the elements~\ref{hof}.
The antisymmetriser $A_k$ is now defined on the $k$-th order component
of $\Gamma^{\otimes}_{\circ}$ which is subsequently factorised by
the kernel of $A_k$ to obtain $\Gamma^{\wedge k}$ from which by direct
summation the differential algebra $\Gamma^{\wedge}$ is formed.
From the properties of $\sigma$ it is seen that in $\Gamma^{\wedge}$
the relations of $\Gamma^{\wedge}_*$ are satisfied.
Thus, $\Gamma^{\wedge}$ is either $\Gamma^{\wedge}_*$ or a factor
algebra of it.

\paragraph{Acknowledgment} The work on this paper has been supported
by a grant of the Studienstiftung des deutschen Volkes.

\emph{Erratum and revision remark:}
In the first version of this paper (as of Feb., 1998), the classification
of freely generated first order differential $*$~calculi was collapsed
by an error in the system of equations used in the proof of the theorem.---%
The second classification result, dealing with differential calculi with
one relation between invariant elements, is new in the revised version.
\end{document}